# QUASI-ESTIMATION AS A BASIS FOR TWO-STAGE SOLVING OF REGRESSION PROBLEM

**Anatoly Gordinsky**
Berman Engineering Ltd
Modiin 71700, Israel

**Abstract**

An effective two-stage method for an estimation of parameters of the linear regression is considered. For this purpose we introduce a certain quasi-estimator that, in contrast to usual estimator, produces two alternative estimates. It is proved that, in comparison to the least squares estimate, one alternative has a significantly smaller quadratic risk, retaining at the same time unbiasedness and consistency. These properties hold true for one-dimensional, multi-dimensional, orthogonal and non-orthogonal problems. Moreover, a Monte-Carlo simulation confirms high robustness of the quasi-estimator to violations of the initial assumptions. Therefore, at the first stage of the estimation we calculate mentioned two alternative estimates. At the second stage we choose the better estimate out of these alternatives. In order to do so we use additional information, among it but not exclusively of a priori nature. In case of two alternatives the volume of such information should be minimal. Furthermore, the additional information is not built-in into the quasi-estimator structure, so that any kind of information, even intuitive one, can be used. These features, in combination with decrease of the quadratic risk, provide a great advantage of our method. A variety of types of the additional information for choosing the better estimate is considered. One example is the successful processing of the famous experiment conducted by astronomers in 1919 to verify the General Theory of Relativity of A. Einstein.

## 1 Introduction

Let us consider the problem of parameters estimation for the following linear regression model:
$$Y = X\beta + \varepsilon, \qquad (1)$$
where $Y$ is a $n \times 1$ vector, $X$ is a $n \times k$ matrix of non-random regressors, $\beta$ is a $k \times 1$ vector of the parameters to be estimated, and $\varepsilon$ is an $n \times 1$ errors vector.

In order to estimate the parameter $\beta$, let us make additional assumptions regarding this parameter, the matrix of regressors $X$ and the vector $\varepsilon$. Suppose the rank of the matrix of regressors is $rank(X) = k$, the errors vector $\varepsilon$ has its mathematical expectation equal to zero, i.e. $E(\varepsilon) = 0$, the variance-covariance matrix of $\varepsilon$ is equal to $cov(\varepsilon) = \sigma^2 I_n$, where $\sigma^2$ is the error variance, $I_n$ - the identity $n$ x $n$-matrix. Regarding the parameter $\beta$ to be estimated, assume that it has no constraints, in other words that $\theta \in R^k$, where $\theta$ - is the set of *a priori* values for $\beta$.

We should make a remark regarding the assumption $E(\varepsilon) = 0$. This assumption can be weakened to $E(\varepsilon) = a$, where $a$ is a constant. In this case, centring $Y$ and $X$ reduces this problem to the one considered. The results of estimating parameters from the centred data are all identical with the results obtained when adding a unit $n \times 1$ column to the matrix $X$. In this case, we have the model (1.1) with an added intercept.

The second remark concerns the requirement $cov(\varepsilon) = \sigma^2 I_n$. If $cov(\varepsilon) = \sigma^2 G$, where $G$ is a known non-identity matrix, then, using a transformation of variables (cf. [1]), we can reduce the model to the initial form.

The above stated conditions completely correspond to the assumptions of the Gauss-Markov theorem (cf. [1]). The theorem states that the least squares estimator (OLS-estimator) defined as

$$b = (X^T X)^{-1} X^T Y, \qquad (2)$$

where "$T$" is the symbol of transposing, is unbiased and has the smallest variance in the class of unbiased estimators linear relative to $Y$.

It is also well-known (cf. [2]) that, for the normal distribution of $\varepsilon$, the estimator (2) will be the best in the class of all (linear and non-linear) unbiased estimators.

It obviously follows from the above that the OLS-estimator, in the mentioned class, has also the smallest mean squared distance from $\beta$, criterion (the quadratic risk) which will be used in the present paper.

However, the fact that the estimator is optimal in one sense does not guarantee that it is good for practical purposes. Let us consider certain situations met in practice. Let us represent the estimator in the following form:

$$b = \beta + \delta, \qquad (3)$$

where $\delta$ is a random variable. Then obtain from (1.2), (1.3) the following equation:

$$\delta = (X^T X)^{-1} X^T \varepsilon, \qquad (4)$$

which, in accordance with properties of $X$ and $\varepsilon$, implies $E(\delta) = 0$.

Let us find now the mean square of the distance from the OLS-estimator to $\beta$, taking into account (3), (4):

$$L^2 = E((b-\beta)^T(b-\beta)) = E(\delta^T \delta) = E(\varepsilon^T X (X^T X)^{-1}(X^T X)^{-1} X^T \varepsilon) =$$
$$= \sigma^2 Sp(X(X^T X)^{-1}(X^T X)^{-1} X^T) = \sigma^2 Sp((X^T X)^{-1}) = \sigma^2 \sum_{i=1}^{k} \lambda_i, \qquad (5)$$

where $Sp$ is the trace of the matrix, $\lambda_i$ is the $i$-th eigenvalue of the inverse matrix $(X^T X)^{-1}$. In obtaining (5) the properties of random quadratic forms and of matrix traces (cf. [1]) were used.

We see from (5) that $L^2$ is increased with the increase in the variance of the errors and the sum of eigenvalues of the inverse matrix, which all have positive values, because the matrix $(X^T X)^{-1}$ is positively defined. The sum of eigenvalues increases with the decrease in number of observations $n$ (the problem of the small sample), and also with the increase in the linear correlation among the columns of the matrix $X$, in other words, in the presence of multicollinearity. The multicollinearity always appears while studying controlled objects, i.e. objects with feedback. Among these there is a vast majority of the operating technological, biological, medical and other objects. The multicollinearity appears also in the absence of feedback, in those cases when the regressors $X$ from the equation (1) are represented by multidimensional polynomials or by some other non-orthogonal series.

There exists also another danger. The structure of the matrix $X$ could be such that in the expression (4) the absolute values of the coefficients of some $\varepsilon_j,...,\varepsilon_m$ will be so big that, even with moderate values of these $\varepsilon$, their impact to the estimate will be too great. Mentioned effects, as well as a number of other effects, can lead to the situation where the estimates $b$ will strongly differ from the true values $\beta$ to be estimated. As a consequence, predicted values of response $\hat{Y}$ calculated by equation $\hat{Y} = Xb$ can be far away from true values $X\beta$.

Precisely because of this, with the aim of increasing quality of the estimation, in the theoretical and applied statistics such methods of the linear regression parameters estimation are being developed and investigated in which, explicitly or implicitly, *a priori* information regarding the parameters to be found are used. Among those are shrinkage estimation, ridge estimation, Bayesian estimation, and others (see e.g. [3], [4], [5]).

In the present paper, the approach, first proposed in [6], [7], [8], is developed in close detail, including new results. This approach is also using a priori or, in general, additional, external information. It differs substantially, however, from the known approaches. The essence of the method is as follows.

The estimation of the regression parameters in the model (1) is conducted in two stages.

At the first stage, the least squares estimation is calculated. Subsequently, it is additively corrected by the vector that has the dimension of vector $\beta$ and depends on a certain random value. This value is the sign of a scalar, determined via the random error (4) of the least squares estimator. The result obtained can be interpreted as a certain quasi-estimator, a function of a discrete random variable taking only two values +1 or -1. With the correct choice of the sign, the quasi-estimate, as is shown further in the paper, will possess a number of useful properties: unbiasedness, consistency, and a significantly smaller value of quadratic risk, compared to the OLS-estimator. Moreover, these properties hold for one-dimensional and multi-dimensional, orthogonal and non-orthogonal problems. The paper also gives the results of the Monte-Carlo simulation, showing high robustness of the quasi-estimator with respect to violations of the initial assumptions. In applications, we have two alternative estimates one of which is superior to the least squares estimate.

At the second stage, the better estimate is chosen out of the two alternatives available. The minimal volume (1 bit) of the required additional information is obvious, and this fact extraordinarily simplifies practical application of the method. Moreover, this information can be of any nature: theoretical, empirical, based on the subject matter experience, obtained as the result of an additional experiment, and even intuitive. It is important that the additional information is not built-in into the estimator construction, and therefore many kinds of information can be used sequentially. All this significantly increases the probability of correct choice and, as a final result, an effectiveness of the method.

A number of types of additional information for choosing the better estimate is reviewed in the paper.

## 2 Quasi-estimator and its properties

Consider two non-linear, non-homogenous in $Y$ estimators:

$$b_1 = b + c\sqrt{e^T e}\, q, \tag{6}$$

$$b_2 = b - c\sqrt{e^T e}\, q, \tag{7}$$

where $b$ is an OLS-estimate (2), $e$ is the ($n \times 1$) vector of the known regression residuals:

$$e = Y - Xb, \tag{8}$$

$q$ is an arbitrary normalized ($q^T q = 1$)  ($k \times 1$) vector, and $c$ is a constant, as yet unknown, that we will subsequently define. We take positive values for the radicals in (6) and (7). The estimates (6) and (7) differ only in the sign of the additive correction to the vector of the OLS-estimator.

*Let us make a convention that, in every application, we choose, out of the two estimators (6) and (7), the one whose correction sign is equal to $sign(-q^T \delta)$, where $\delta$ is the random error for the OLS-estimator, determined by (4), and the sign-function is defined by the rule: $sign(x) = 1$ if $x \geq 0$, and $sign(x) = -1$ if $x < 0$.*

Now we have obtained the following quasi-estimator:

$$\tilde{b} = b - sign(q^T \delta) c \sqrt{e^T e}\, q \qquad (9)$$

The term "quasi-estimator" is used because the expression (9) involves a discrete random variable $sign(q^T \delta)$, which receives just two values: +1 or -1.

Let us determine the constant $c$ in such a way so that the average square of the distance between the quasi-estimator (9) and $\beta$ is minimized for any $q$. The aforementioned mean square of the distance, for $\tilde{b}$, is equal to:

$$\tilde{L}^2 = E((\tilde{b} - \beta)^T (\tilde{b} - \beta)). \qquad (10)$$

Substituting into (10) the value of $\tilde{b}$ from (9), and then of $b$ from (3), taking the derivative with respect to $c$, and equating the result to zero, we obtain the following:

$$\tilde{c} = \underset{c \in R^1}{arg(min\, \tilde{L}^2)} = E(\sqrt{(q^T \delta)^2 e^T e}) / E(e^T e). \qquad (11)$$

Now the quasi-estimate with the minimal $\tilde{L}^2$ will have the following form:

$$\tilde{b} = b - sign(q^T \delta) \cdot \tilde{c} \sqrt{e^T e} \cdot q, \qquad (12)$$

where $\tilde{c}$ is taken from (11).

Note that when $n = k$ the residual sum of squares is equal to $e^T e = 0$ and, as follows from (12), $\tilde{b} = b$.

Thus, in the sequel, we will restrict ourselves to the case $n > k$.

The aforementioned is already sufficient to prove the following:

**PROPOSITION 1.** *Suppose $E(\varepsilon) = 0$, $cov(\varepsilon) = \sigma^2 I_n$, $n > k$ and $q$ is an arbitrary normalized $k \times 1$ vector. Then $\tilde{L}^2 < L^2$, for the quasi-estimator (12).*

**Proof.** Substitute (11) into (12), and the result obtained – into (10). After some calculations, writing $sign(q^T \delta) \cdot q^T \delta = |q^T \delta| = \sqrt{(q^T \delta)^2}$ and using (3) and (5), we obtain:

$$\tilde{L}^2 = L^2 - \left[ E(\sqrt{(q^T \delta)^2 e^T e}) \right]^2 / E(e^T e) \qquad (13)$$

Since the value being subtracted from $\tilde{L}^2$ is positive, the proposition is proven. ∎

Up to this point, we did not consider the distribution of the error $\varepsilon$. Assume now that $\varepsilon$ is normally distributed, i.e. that $\varepsilon \in N(0, \sigma)$, maintaining all our previous stipulations. In that case we can, first of all, determine what the vector $q$ should be in order to minimize the value $\tilde{L}^2$ from (13).

**PROPOSITION 2.** *Suppose $E(\varepsilon) = 0$, $cov(\varepsilon) = \sigma^2 I_n$, $n > k$ and the distribution of the error $\varepsilon$ is normal, $\varepsilon \in N(0, \sigma)$.*

*Then the minimum of $\tilde{L}^2$ is achieved when the vector $q$ from (12) equals the normalized eigenvector $z_1$ ($z_1^T z_1 = 1$) which corresponds to the maximum eigenvalue of the inverse matrix $(X^T X)^{-1}$:*

$$z_1 = \arg(\min_{\substack{q \in R^k \\ q^T q = 1}} \tilde{L}^2)$$

**Proof**

Let us consider (13) and, first, show independence of the quadratic forms $(q^T \delta)^2$ and $e^T e$. Toward this, using (2), (4), (8) and, taking into account that $q^T \delta$ is a scalar, let us represent these forms as $\varepsilon^T T \varepsilon$ and $\varepsilon^T B \varepsilon$, correspondingly, where the matrices $T$ and $B$ are equal, correspondingly, to:

$$T = X(X^T X)^{-1} q q^T (X^T X)^{-1} X^T, \qquad (14)$$

$$B = I_n - X(X^T X)^{-1} X^T \qquad (15)$$

Direct verification shows that $TB = BT = 0$. This is a necessary and sufficient condition for independence of the quadratic forms from the considered normal random values $\varepsilon$ (cf. [9]).

It follows from the independence, that the denominator in (13) (while the radicals are taken positive) can be represented as a product:

$$\left[E(\sqrt{(q^T \delta)^2 e^T e})\right]^2 = \left[E(\sqrt{(q^T \delta)^2})\right]^2 \left[E(\sqrt{e^T e})\right]^2 = \left[E(|q^T \delta|)\right]^2 \left[E(\sqrt{e^T e})\right]^2.$$

In order to minimize $\tilde{L}^2$, one needs to maximize the first multiplicand of the above product. Scalar $q^T \delta$ is a normal random variable with the zero mathematical expectation.

The module of such a variable has semi-normal distribution. Its mathematical expectation is known (cf. [10]):

$$E(|q^T \delta|) = \sqrt{\frac{2}{\pi} D(q^T \delta)} \Rightarrow \left[E(|q^T \delta|)\right]^2 = \frac{2}{\pi} D(q^T \delta),$$

where $D$ is the symbol of variance. Let us find this variance:

$$D(q^T \delta) = E(\delta^T q q^T \delta) = E(\varepsilon^T X(X^T X)^{-1} q q^T (X^T X)^{-1} X^T \varepsilon) =$$
$$= \sigma^2 Sp(X(X^T X)^{-1} q q^T (X^T X)^{-1} X^T) = \sigma^2 (q^T (X^T X)^{-1} q),$$

where, as designated above, $Sp$ is the symbol of the trace of the matrix and $\sigma$ is the standard deviation of $\varepsilon$. Obtaining this result, we use properties of the matrix traces and the equality of the trace of the scalar to the scalar itself.

To finish the proof of the proposition, we are left to find out what should be the value of the normalized vector $q$ in order to maximize the value of the obtained scalar. The answer is found right in the Rayleigh-Ritz theorem (cf. [11]), according to which this maximum is equal to the maximal eigenvalue $\lambda_1$ of the matrix $(X^T X)^{-1}$ and is obtained at the normalized eigenvector of this matrix corresponding to such a maximum eigenvalue.
The proposition is proven. ∎

Now the quasi-estimator gets the following form:

$$\tilde{b} = b - sign(z_1^T \delta) \cdot \tilde{c} \sqrt{e^T e} \cdot z_1 \quad (16)$$

Normality of $\varepsilon$ and the above result allow us to define concretely the expression $\tilde{c}$ in (16) and, eventually, obtain the final form for the optimal quasi-estimator.

**PROPOSITION 3.** *Let* $E(\varepsilon) = 0$, $cov(\varepsilon) = \sigma^2 I_n$, $n > k$ *and* $\varepsilon \in N(0,\sigma)$.
*Then the quasi-estimator (16) gains the following form:*

$$\tilde{b}_o = b - sign(z_1^T \delta) \sqrt{\frac{\lambda_1}{\pi}} \frac{\Gamma((n-k+1)/2)}{\Gamma((n-k+2)/2)} \sqrt{e^T e}\, z_1, \quad (17)$$

where $\Gamma(x)$ is the gamma-function, and $\lambda_1$ and $z_1$, as noted above, are the maximum eigenvalue of the matrix $(X^T X)^{-1}$ and the normalized eigenvector corresponding to it, respectively.

**Proof.**
Let us find $\tilde{c}$ from (11) and substitute the result into (16).
First let us find the numerator in (11). Taking into account (4) and the property $(X^T X)^{-1} z_1 = \lambda_1 z_1$ of the eigenvector of the matrix (cf. [1]), let us represent $(z_1^T \delta)^2$ in the following form: $(z_1^T \delta)^2 = \sigma^2 \lambda_1 \tilde{\varepsilon}^T A \tilde{\varepsilon}$, where $\tilde{\varepsilon}$ is the normalized error, $\tilde{\varepsilon} = \varepsilon / \sigma$, $\tilde{\varepsilon} \in N(0,I)$ and $A = \lambda_1 X z_1 z_1^T X^T$. Checking shows that $A^2 = A$, i.e. $A$ is an idempotent matrix with $rank(A) = Sp(A) = 1$. When $\tilde{\varepsilon} \in N(0,I)$, this is a necessary and sufficient condition for the assertion $\tilde{\varepsilon}^T A \tilde{\varepsilon} \in \chi_1^2$, i.e. the shown quadratic form has distribution $\chi^2$ with one degree of freedom. Similarly represent the quadratic form $e^T e = \sigma^2 \tilde{\varepsilon}^T B \tilde{\varepsilon}$, where the matrix $B$ is determined from (15).
Then $B = I_n - X(X^T X)^{-1} X^T$, $B^2 = B$ and $rank(B) = n - k$, whence $\tilde{\varepsilon}^T B \tilde{\varepsilon} \in \chi_{n-k}^2$.
The numerator in the expression for $\tilde{c}$ will get now the following form:

$$E(\sqrt{(z_1^T \delta)^2 e^T e}) = \sigma^2 \sqrt{\lambda_1}\, E(\sqrt{\tilde{\varepsilon}^T A \tilde{\varepsilon} \cdot \tilde{\varepsilon}^T B \tilde{\varepsilon}}).$$

Analogously to our earlier remark, the quadratic forms under the sign of radical are independent, since $AB = BA = 0$.

Let us now compute the mathematical expectation of the positive value for the square root of the product of independent random variables distributed as $\chi_1^2$ and $\chi_{n-k}^2$. Using known relations between the distribution densities and the mathematical expectation of the function (cf. [10]), we obtain:

$$E(\sqrt{(z_1^T \delta)^2 e^T e}) =$$

$$= \sigma^2 \frac{\sqrt{\lambda_1}}{\sqrt{2\pi}} \int_0^\infty x^{0.5} x^{-0.5} \exp(-\frac{x}{2}) dx \frac{1}{2^{(n-k)/2} \Gamma((n-k)/2)} \int_0^\infty y^{0.5} y^{(n-k-2)/2} \exp(-\frac{y}{2}) dy =$$

$$= \sigma^2 \frac{\sqrt{\lambda_1}}{\sqrt{2\pi}} 2 \frac{\Gamma((n-k+1)/2)}{2^{(n-k)/2} \Gamma((n-k)/2)} 2^{(n-k+1)/2} = \frac{2\sigma^2 \sqrt{\lambda_1}}{\sqrt{\pi}} \frac{\Gamma((n-k+1)/2)}{\Gamma((n-k)/2)}.$$

Substituting the mathematical expectation $e^T e$ into the denominator in the expression (11), namely $E(e^T e) = (n-k)\sigma^2$, and using the property $\Gamma(x+1) = x\Gamma(x)$. (cf. [11]) of the gamma-function, we obtain $\Gamma((n-k+2)/2) = ((n-k)/2)\Gamma((n-k)/2)$ and

$$\tilde{c} = \frac{\sqrt{\lambda_1}}{\sqrt{\pi}} \frac{\Gamma((n-k+1)/2)}{\Gamma((n-k+2)/2)} \qquad (18)$$

We are left to substitute this expression for $\tilde{c}$ in the equation (16), and Proposition 3 has been proven. ∎

Possessing the definitions and the results (6), (11), (18), we can compute the ratio of the $\tilde{L}_o^2$ for the optimal quasi-estimator to the $L^2$ of the OLS-estimator:

$$\tilde{L}_o^2 / L^2 = 1 - \frac{n-k}{\pi} \frac{\Gamma^2((n-k+1)/2)}{\Gamma^2((n-k+2)/2)} \frac{\lambda_1}{\sum_{i=1}^{k} \lambda_i} \qquad (19)$$

Using the relation from [11]:

$$\Gamma(\alpha + p)/\Gamma(\alpha + h) = \alpha^{p-h}(1 + \frac{1}{2\alpha}(p-h)(p+h-1)) + O(1/\alpha^2),$$

we can get from (19) the equation that is more easy for the following analysis when $(n-k)^2/4 \gg 1$:

$$\tilde{L}_o^2 / L^2 = 1 - \frac{2}{\pi}(1 - \frac{0.25}{n-k})^2 \frac{\lambda_1}{\sum_{i=1}^{k} \lambda_i} \qquad (20)$$

As can be seen from (19) and (20), the relative gain for the quasi-estimator depends mainly on the distribution of the eigenvalues $\lambda_i$ of the matrix $(X^T X)^{-1}$, and can be quite substantial. Thus, in the one-dimensional case, when $\lambda_1 = \sum_{i=1}^{k} \lambda_i$, the ratio $\tilde{L}_o^2 / L^2$ is close to 0.4, i.e. the quasi-estimator has $\tilde{L}_o^2$ smaller than the OLS-estimator $L^2$ by the ratio of 2.5. This fact can have a great significance, in particular, when processing direct measurements.

In the case of multicollinearity, when $\lambda_1 \cong \sum_{i=1}^{k} \lambda_i$, we obtain approximately the same result.

When processing orthogonal data, the effect, naturally, will be smaller and will substantially depend on the number of variables $k$.

Let us establish some more important properties of the optimal quasi-estimator. When the property under consideration holds for the original quasi-estimator (9), with an arbitrary vector $q$, we will mention this in our remarks.

**PROPOSITION 4** *Let $E(\varepsilon) = 0$, $cov(\varepsilon) = \sigma^2 I_n$, $n > k$ and $\varepsilon \in N(0,\sigma)$. Then the optimal quasi-estimator $\tilde{b}_o$ from (17) is unbiased, i.e. $E(\tilde{b}_o) = \beta$.*

**To prove this**, let us compute the mathematical expectations of both sides of the equation (17), taking into account the independence of $z_1^T \delta$ and $e^T e$ proved above and therefore of their functions:

$$E(\tilde{b}_O) = E(b) - E(sign(z_1^T \delta))E(\sqrt{e^T e}) \sqrt{\frac{\lambda_1}{\pi}} \frac{\Gamma((n-k+1)/2)}{\Gamma((n-k+1)/2)} z_1.$$

Using (4) and the defining property of the eigenvector, we obtain:
$$z_1^T \delta \in N(0, \sqrt{\lambda_1}\,\sigma), \tag{21}$$

that is, the random variable $z_1^T \delta$ has normal distribution with the zero mathematical expectation and variance $\lambda_1 \sigma^2$. The mathematical expectation of this function is equal to the following:

$$E(sign(z_1^T \delta)) = \int_{-\infty}^{\infty} sign(x) f(x)\, dx, \text{ where } f(x) \text{ - is the distribution density of the normally}$$

distributed value (21), continuous, symmetrical and bounded above.
In view of our definition of the *sign*-function ( $sign(x) = 1$ when $x \geq 0$ and $sign(x) = -1$ when $x < 0$ ), the given integral is equal to zero due to the symmetry of the integrand function for all $x \neq 0$ and its boundedness at the point $x = 0$.

From this, since $E(sign(z_1^T \delta)) = 0$ and, because of unbiasedness of the OLS-estimator $E(b) = \beta$, we get $E(\tilde{b}_O) = \beta$. ∎

*Remark 1*. We will obtain the same result for any continuous, symmetric, and bounded above distribution of the error $\varepsilon$, for which the mathematical expectation exists.

*Remark 2*. The property of unbiasedness, under conditions of Remark 1, applies also to the estimator (9).

**PROPOSITION 5.** *Let* $E(\varepsilon) = 0$, $cov(\varepsilon) = \sigma^2 I_n$, $n > k$ *and* $\varepsilon \in N(0, \sigma)$.

*Then the variance-covariance matrix for the optimal quasi-estimator (17) is equal to:*

$$Q = E((\tilde{b}_O - E(\tilde{b}_O))(\tilde{b}_O - E(\tilde{b}_O))^T) = \sigma^2 \left[ (X^T X)^{-1} - (n-k)\frac{\lambda_1}{\pi} \frac{\Gamma^2((n-k+1)/2)}{\Gamma^2((n-k+2)/2)} z_1 z_1^T \right] \tag{22}$$

**Proof.**
Proceeding from the definition of the variance-covariance matrix given on the left side of (22), and taking into account the unbiasedness of the estimator $\tilde{b}_O$, and also the relations (3),(17),(18), after some manipulations, we obtain:

$$Q = E(\delta\delta^T - \tilde{c}\, sign(z_1^T \delta) z_1 \delta^T \sqrt{\varepsilon^T B \varepsilon} - \tilde{c}\, sign(z_1^T \delta) \delta z_1^T \sqrt{\varepsilon^T B \varepsilon} + \tilde{c}^2 \varepsilon^T B \varepsilon\, z_1 z_1^T) \tag{23}$$

Let us consider the second summand in (23) which includes the product of three random variables.

Using (4) and the known expansion $(X^T X)^{-1} = \sum_{i=1}^{k} \lambda_i z_i z_i^T$, express this summand as follows:

$$ad2 = -\tilde{c}\, sign(z_1^T \delta)\, z_1 \varepsilon^T X \sum_{1}^{k} \lambda_i z_i z_i^T \sqrt{\varepsilon^T B \varepsilon}. \tag{24}$$

Let us prove that when $i > 1$ the mathematical expectation of the corresponding summand in (24) is equal to zero. Let us first establish the mutual independence of the three random scalars $z_1^T \delta$, $\varepsilon^T X z_i$, $\sqrt{\varepsilon^T B \varepsilon}$ and, therefore, their functions. For this purpose use the following fact (cf. [9]): *for positive semi-definite quadratic forms of normal random variables, their statistical independence follows from their pairwise independence.* The independence of the first and the third scalars was established during the proof of Proposition 2. The independence of the first and the second normally distributed scalars follows from their noncorrelatedness: $E(z_1^T X^T \varepsilon \varepsilon^T X z_i) = \sigma^2 z_1^T X^T X z_i = \sigma^2 \lambda_i z_1^T z_i = 0$ because of the orthogonality of the eigenvectors. Finally, let us prove the independence of the second scalar and the quadratic form $\varepsilon^T B \varepsilon$. For this purpose, divide the scalar by $\sigma$, multiply by $\sqrt{\lambda_i}$ and form from the result a quadratic form $\tilde{\varepsilon}^T A \tilde{\varepsilon}$, where $\tilde{\varepsilon} \in N(0, I)$ and $A = \lambda_i X z_i z_i^T X^T$. Verification shows that the matrix $A$ is idempotent (the idempotency of $B$ has been established earlier, during the proof of Proposition 3) and that $AB = 0$, from which follows independence of the second and the third scalars in (24), and the mutual independence of the members of the triple product as a whole, when $i > 1$. In this case, the mathematical expectation of the product is equal to the product of the mathematical expectations, and, as two of those are equal to zero, it is also equal to zero.

Let us find the mathematical expectation of the first summand in $ad2$ (24), corresponding to $i = 1$:

$$-\tilde{c} E( sign(z_1^T \delta) z_1 \varepsilon^T X \lambda_1 z_1 z_1^T \sqrt{\varepsilon^T B \varepsilon} ) = -\tilde{c} E( sign(z_1^T \delta) z_1^T \delta \sqrt{\varepsilon^T B \varepsilon} ) z_1 z_1^T =$$
$$= -\tilde{c} E(\sqrt{(z_1^T \delta)^2 \varepsilon^T B \varepsilon}) z_1 z_1^T.$$

The obtained mathematical expectation has been found during the proof of Proposition 3. Using this result, the expression for the gamma-function used in the same proof, and the expression (18), we obtain: $ad2 = -\tilde{c}^2 \sigma^2 (n-k) z_1 z_1^T$. Using the same method, we get the same result for the third summand in (23). Finally, the mathematical expectation of the first and fourth summands in (23) is known:

$$E(\delta \delta^T + \tilde{c}^2 \varepsilon^T B \varepsilon z_1 z_1^T) = \sigma^2 ((X'X)^{-1} + \tilde{c}^2 (n-k) z_1 z_1^T).$$

Putting the obtained results together and substituting $\tilde{c}$ from (18), we will get (22). ∎

*Remark.* It follows from (22) that, under the above mentioned assumptions, every element of the vector of the quasi-estimator $\tilde{b}_o$ has variance not greater than the variance of the corresponding element of the vector $b$ of the OLS-estimator. That is, taking into account the unbiasedness of the quasi-estimator, we can consider this result as a certain analogue, for the optimal quasi-estimator, of the Gauss-Markov's theorem, however, in a restricted form, not for an arbitrary distribution, but only for normal distribution of the error $\varepsilon$.

Let us note right away that, in the one-dimensional case, as it follows from Proposition 1, the quasi-estimator (even in its original form, and still more for the optimal quasi-estimator) always has smaller variance than the OLS-estimator, whatever the distribution of the error is.

Let us present also the following facts regarding the variance-covariance matrix $Q$. It is easy to verify that the maximal eigenvalue of the matrix $Q$ is equal to:

$$\lambda_Q = \lambda_1(1 - \frac{n-k}{\pi}\frac{\Gamma^2((n-k+1)/2)}{\Gamma^2((n-k+2)/2)}), \qquad (25)$$

that all other eigenvalues are equal to the corresponding eigenvalues of the matrix $(X^T X)^{-1}$, all eigenvectors of the matrix $Q$ are equal to the eigenvectors of the matrix $(X^T X)^{-1}$. In addition, obviously, the matrix $Q$ is positively defined.

Let us establish yet another property of the optimal quasi-estimator (17).

**PROPOSITION 6.** *If the OLS-estimator is consistent in the mean-square sense, then the optimal quasi-estimator (17) is also consistent in the same sense.*

**Proof.** Suppose that $\lim_{n\to\infty}(X^T X)_n^{-1} = 0$, i.e. that the OLS-estimator is quadratic mean consistent (cf. [9]). Consider the variance-covariance matrix (22) for the optimal quasi-estimator $\tilde{b}_O$ as a function of $n$.

Substitute $(X^T X)_n^{-1} z_{1n}$ instead of $\lambda_{1n} z_{1n}$ and take out $(X^T X)_n^{-1}$:

$$Q_n = \sigma^2 (X^T X)_n^{-1}(I_k - \frac{n-k}{\pi}\frac{\Gamma^2((n-k+1)/2)}{\Gamma^2((n-k+2)/2)} z_{1n} z_{1n}^T).$$

Using the relation for the gamma-function first used in the derivation of (20), we obtain:

$$\lim_{n\to\infty}((n-k)\frac{\Gamma^2((n-k+1)/2)}{\Gamma^2((n-k+2)/2)}) = 2.$$

Taking into account that $z_{1n} z_{1n}^T$ is bounded, since the eigenvector $z_{1n}$ is normalized, we finally get $\lim_{n\to\infty} Q_n = 0$. ∎

Let us now consider questions related to estimation of the confidence intervals for the optimal quasi-estimator (17). Taking into account its unbiasedness, represent the estimator in the following form: $\tilde{b}_o = \beta + \tilde{\delta}$, where, in view of (3), (17), (18)

$$\tilde{\delta} = \delta - sign(z_1^T \delta)\,\tilde{c}\, z_1 \sqrt{e^T e} \qquad (26)$$

Let us consider the moments of the *j*-component of the vector $\tilde{\delta}$, denoting it by $\tilde{\delta}_j$. Remember that the mathematical expectation of the component is 0. Short of producing a tedious proof, let us only remark that the odd central moments of this component are equal to 0, and the even ones are not much different from the normal distribution moments. In particular, the variance of this distribution is equal to
the j-th diagonal element of the variance-covariance matrix $Q$ from (22), and the forth momentum is determined by the following equation:

$$\frac{\mu_4(j)}{\sigma^4} = z_1(j)^4\left[3\lambda_1^2 - 2(n-k)\lambda_1\tilde{c}^2 - (n-k)(3(n-k)+2)\tilde{c}^4\right] +$$
$$+ 6z_1(j)^2\left[\lambda_1 - (n-k)\tilde{c}^2\right]A_j^T A_j + 3(A_j^T A_j)^2,$$

where $\quad A_j = \sum_{i=2}^{k}\lambda_i z_i(j) X z_i$, and also $A_j = 0$, if $k=1$. $\qquad (27)$

Dividing $\mu_4(j)$ from (27) by the squares of the corresponding diagonal elements of the matrix $Q$ and subtracting 3, we are obtaining the kurtosis of the component $\tilde{\delta}_j$. The analysis shows that the excess is positive, varying in the interval between 0 and 1, depending on the characteristics of the matrix $X$, and reaching the upper limit at the components with a high co-variance. This means that the distribution is somewhat more leptokurtic than the normal one. Therefore, when assuming it to be approximately normal, we will have a reserve in the confidence interval.

Thus, using the approximate relation

$$\tilde{b}_o \in N(\beta, Q), \qquad (28)$$

we can extend all known results from the theory of the regressive analysis (cf. [2], [9]) to the optimal quasi-estimator (17). In particular, the individual confidence interval for the $j$-component of the vector of the quasi-estimator is computed with the help of the formula:

$$\tilde{b}_o(j) \pm t(n-k, 1-\alpha/2)\sqrt{Q_{j,j}}\, s, \qquad (29)$$

where $t(n-k, 1-\alpha/2)$ is the $1-\alpha/2$ point of the Student distribution with $n-k$ degrees of freedom, $Q_{j,j}$ - is the corresponding diagonal element of the matrix $\dfrac{Q}{\sigma^2}$, $s$ - is the standard deviation estimate:

$$s = \sqrt{e^T e/(n-k)} \qquad (30)$$

The joint $100(1-\alpha)$-percent confidence area for the quasi-estimator $\tilde{b}_o$ will be determined from the expression:

$$(\beta - \tilde{b}_o)^T Q (\beta - \tilde{b}_o) \leq k s^2 F(k, n-k, 1-\alpha), \qquad (31)$$

where $F(k, n-k, 1-\alpha)$ is the $(1-\alpha)$-percent point for the $F(k, n-k)$ distribution.

Finishing consideration of the basic properties of the optimal quasi-estimator (17), let us answer the following important question: how the efficiency index (19) will change when the original assumptions are violated, namely if the error distribution $\varepsilon$ is symmetrical, but different from normal, and the error covariance matrix is not the identity matrix.

An answer to this question was obtained with the help of statistical modelling using Statistics Toolbox from the Matlab package. Errors with the normal distribution $N(0,1)$, uniformly distributed in the interval (-2,2), a mixture of the normal distributions $N(0,1)$-80% and $N(0,10)$-20%, and also the errors representing a time series with the exponential autocorrelation function $R(\tau) = \sigma^2 e^{-q|\tau|}$, with $\sigma^2 = 1$ and $q = 0.3$, were modelled. The number of tests was $10,000$. As data, the unit vector $X_{10,1}$ has been used (i.e. a one-dimensional problem of direct measurements $Y = \beta + \varepsilon$ was modelled), as well as the matrix $X_{10,2}$ with a significant linear contingency of its columns (i.e. a two-dimensional problem $Y = X_1 \beta_1 + X_2 \beta_2 + \varepsilon$ was modelled under conditions of multicollinearity). The degree of multicollinearity is characterized by the following two indicators: the ratio of the maximal eigenvalue of the matrix $\left(X^T X\right)^{-1}$ to the minimal one, which was equal to 448.8, and the correlation coefficient between the columns of the matrix $X$, which was equal to 0.9955. The data for the two-dimensional problem were centred. The results of modelling are shown in the tables 1 and 2 below.

*Table 1. One-dimensional problem*

|  | Normal distribution N(0.1) | Uniform distribution in interval (-2÷2) | Mixture of distributions N(0,1)-80% and N(0,10)-20% | Exponential autocorrelation $\sigma = 1, q = 0.3$ |
|---|---|---|---|---|
| *Theoretical value of $\tilde{L}^2 / L^2$* | 0.39772 | - | - | - |
| *Experimental value of $\tilde{L}^2 / L^2$* | 0.39472 | 0.41411 | 0.25095 | 0.65278 |

*Table 2. Two-dimensional problem*

|  | Normal distribution N(0.1) | Uniform distribution in interval (-2÷2) | Mixture of distributions N(0,1)-80% and N(0,10)-20% | Exponential autocorrelation $\sigma = 1, q = 0.3$ |
|---|---|---|---|---|
| *Theoretical value of $\tilde{L}^2 / L^2$* | 0.40319 | - | - | - |
| *Experimental value of $\tilde{L}^2 / L^2$* | 0.40013 | 0.40937 | 0.37667 | 0.52594 |

One can see from Tables 1 and 2 that the optimal quasi-estimator shows a high degree of robustness, i.e. of maintaining stability in the presence of deviations from original assumptions regarding the error $\varepsilon$. A visible reduction in efficiency is observed only in the one-dimensional problem with the autocorrelated error $\varepsilon$, which is explained by the absence of centred data in the one-dimensional case. On the other hand, in the one-dimensional case, and with very heavy characteristics of the mixture of the distributions, one observes an increase in efficiency, i.e. in this case the robustness of the quasi-estimator is increased. For the two-dimensional problem, almost the same picture is typical, but with smaller scatter for the efficiency criterion $\tilde{L}^2 / L^2$.

## 3 Choice of a better estimate from two alternatives

Let us next consider some of the possibilities for the choice of the best estimates from two alternative ones. For this purpose, the additional information available to the researcher is used. This information can be of any nature: theoretical, empirical, based on the subject matter experience, obtained from an additional experiment of a smaller volume, based on one's intuition.

Here there are two reasons allowing us to assert that we are obtaining maximal results from using additional information. First, other conditions being equal, we need minimal additional information (1 bit) for the choice of one out of the two estimates; second, we are able to use all kinds of additional information, available to the researcher, sequentially, because it is not incorporated into the structure of the quasi-estimator. And then, in different applications, different kinds of *a priori* information can be efficient.

3.1 *Case study: confirmation of one of the two competing theories with the help of an experiment.*

As an example, let us consider a widely known experiment conducted by astronomers Dyson, F.W., Eddington, A.S., and Davidson, C.R. in the year of 1919 (cf. [13]). The purpose of the experiment was to determine the deflection of a ray of light in the gravitational field of the sun. They considered three possibilities: the deflection is absent; the deflection conforms to Newton's theory and is equal to 0.87"; the deflection conforms to the General Theory of Relativity of A. Einstein and is equal to 1.75". As the outcome of the experiment, there were obtained three independent values: $1.98 \pm 0.12$, $1.61 \pm 0.3$, and 0.93. For the first two values, after the symbol $\pm$, the probable error is shown, which is equal, as is known, to $0.6745\sigma$. For the third value, the authors did not give the value of the probable error. Noting that it is too big, they discarded that value. Subsequently, this fact led to doubts and prolonged discussions, and only a double checking, conducted in 1979, has shown that the error in the third measurement was indeed big (cf. [14]). We will process these data without discarding the third measurement 0.93 and assuming its probable error as equal to $\pm 0.6$. Taking into account independence and using a known transformation (cf. [15]) to remove their heteroscedasticity, we obtain, according to (17), the alternative estimates: 2.0051" and 1.7862". The choice of one of these estimates is obvious, since the estimate 2.0051 exceeds all *a priori* allowed values. Thus the final result is 1.7862", which is very close to the value 1.75", given by the General Theory of Relativity.

It is important to note that, even if we substantially reduce the expected probable error for the third measurement, changing it to $\pm 0.4$, i.e. assigning greater weight to the measurement 0.93, we will get the following alternative estimates: 2.0 and 1.7141. The final estimate, in this case, would be 1.7141, which is also close to the theoretical value. Let us also note that the last fact demonstrates the robustness of the quasi-estimator.

3.2 *Case study: studying systems with negative feed-back and known parameters' signs.*

Suppose that we are to estimate the increase of blood glucose level in a patient with the insulin-dependent form of *diabetes mellitus* when he/she is receiving a certain quantity of carbohydrates orally. In particular, precisely such a situation arises when conducting tests on glucose tolerance. However, such a test induces the patient into a state of temporary, and sometimes prolonged, hyperglycemia and, therefore, of decompensation, which is extremely undesirable. Avoiding such a situation is possible by injecting the patient with insulin at the same time as glucose is taken in.

However, in this case, the statistician processing these data will encounter two difficulties. First, the change in the blood sugar level will be insignificant; therefore the ratio signal/noise will be big. Second, the prediction variables will be strongly correlated, i.e. the statistician will encounter the problem of multicollinearity.

Such a situation is typical when studying all systems with the negative feedback, which means all functioning biological, technological and others systems. Precisely because of this, for example, it is very complex and, even impossible, to conduct active experiments in industry according to the plans developed in the theory of experimental design. These experiments are destroying the functioning and necessary negative feedbacks, which leads, in the best case, to the production of defective goods, and in the worse case, to accidents, including the grave ones (cf., for example, [16]).

In the considered situation, application of the quasi-estimator proves to be quite effective. In addition, the choice of the better estimate, out of the two, becomes obvious, since the signs of the parameters are known. For example, in the case of studying parameters of the diabetic patient, it is obvious that a greater amount of the orally received glucose leads to higher blood glucose levels, which means that the corresponding parameter is positive, while, with the increase in the quantity of insulin, the blood glucose level is decreased, and therefore the corresponding parameter is negative.

As an illustration for a study of a diabetes patient, consider a simplified example of the form $Y = X_1 * \beta_1 + X_2 * \beta_2 + \varepsilon$, where $X_1$ - is the increment in the glucose level obtained, relative to the mean of the sample, $X_2$ - the analogues increment for the quantity of insulin, $Y$ - the same relative to the mean for the blood glucose, $\varepsilon$ - is the random error which is an application of the normal distribution $N(0,1)$, also centred relative to the mean. Suppose now that there exists someone who knows the patient's parameters and suppose that they are equal to $\beta_1$ = **40**, $\beta_2$ = **-37,** which means that the vector $\beta$ is equal to $\beta = \begin{bmatrix} 40 \\ -37 \end{bmatrix}$. Also suppose that now this someone wants to see what estimates will be obtained by us applying our approach. Exactly this kind of method of testing of the quality of estimation is the strictest. The data for the given problem is collected in Table 3.

*Table 3. Example of a two-dimensional regression $\beta_1$ = 40, $\beta_2$ = -37*

| # | 1 | 2 | 3 | 4 | 5 | 6 | 7 | 8 | 9 | 10 |
|---|---|---|---|---|---|---|---|---|---|---|
| $X_1$ | -0.061 | -0.051 | 0.059 | -0.271 | 0.109 | 0.099 | -0.101 | 0.149 | 0.089 | -0.021 |
| $X_2$ | -0.055 | -0.055 | 0.065 | -0.255 | 0.085 | 0.095 | -0.085 | 0.155 | 0.065 | -0.015 |
| $\varepsilon$ | 0.3132 | 0.9672 | 1.5252 | -0.7748 | -1.0008 | -1.6578 | -0.4138 | 1.6742 | -0.3428 | -0.2898 |
| $Y$ | -0.0918 | 0.9622 | 1.4802 | -2.1798 | 0.2142 | -1.2128 | -1.3088 | 1.8992 | 0.8122 | -0.5748 |

It is easy to see that $X_1$ and $X_2$ are strongly correlated, i.e. the problem is characterised by multicollinearity. Let us now perform a standard regression analysis and obtain the OLS-estimate. $b = \begin{bmatrix} -0.2868 \\ 7.9614 \end{bmatrix}$. Note that the received estimate $b$ drastically differs from true values $\beta$, and the norm of the obtained OLS-estimations is sharply smaller than the norm of the true coefficients: $b^T b << \beta^T \beta = 2969$. This is wholly typical for systems of this kind, those with negative feedback, in spite of the well-known (and also easily obtainable from the preceding) relation $E(b^T b) = \beta^T \beta + \sigma^2 \sum_{i=1}^{k} \lambda_i > \beta^T \beta$. Thus, for the data collected in the Table 3, if one starts to vary $\varepsilon$, the norm of the OLS-estimations will become smaller than the norm of the true coefficients, approximately in 50% of all cases. It is also obvious that the OLS-estimations are remote from the true ones and that the square of the Euclidian distance is $(b - \beta)^T (b - \beta) = 3644.6$.

Now, we obtain, using the relation (17), two alternative results corresponding to the possible values $sign(z_1^T \delta) = +1$ or -1: $b_1 = \begin{bmatrix} 21.473 \\ -15.329 \end{bmatrix}$ and $b_2 = \begin{bmatrix} -22.047 \\ 31.252 \end{bmatrix}$.

However, the signs of the coefficients are known to us, and therefore it is natural to choose $b_1$ as the correct possibility. It is easy to verify that the square of the Euclidean distance to the true values for the obtained estimate will be equal to 812.9, in other words, 4.48 times smaller.

3.3 *Case study: known constrains on the parameters, given in the form of inequalities.*
In order to formalize the use of additional information, let us consider, from all of types of information and all approaches corresponding to them, the approach based on the availability of the additional information in the area where the unknown parameters are lying. Additionally, let us restrict ourselves to the one-dimensional case. Such a choice is justified by the fact that the results obtained can be extended to the multi-dimensional case, because a multi-dimensional task can be reduced to a certain set of one-dimensional tasks.

When using the mentioned additional information, the final estimator $\overset{\bullet}{b}$ can possess properties different from the properties of the quasi-estimator $\tilde{b}_o$. Thus $\overset{\bullet}{b}$ can become a biased estimate with its variance – higher than the variance of $\tilde{b}_o$. However, the mean sum of the squares of the deviations of $\overset{\bullet}{b}$ will be lower than the same for the OLS-estimation. However, it is reasonable to compare properties of the resulting estimate $\overset{\bullet}{b}$ with the properties of that OLS-estimation which uses the same a priori information, in other words, with the OLS-estimation with the constrains in the form of inequalities. As will be shown below, compared to this one, the proposed estimate $\overset{\bullet}{b}$ will also have an advantage.

Let $Y = X_1 \beta_1 + \varepsilon$, *where* $X_1$ - is a $n \times 1$ column vector. As is known, the variance of the OLS-estimation $b$ is in this case $\tilde{\sigma}^2 = \sigma^2 / \sum_{i=1}^{n} x_{i,1}^2$. Suppose that it is known that $a_1 \leq \beta_1 \leq a_2$. Minimizing the sum of the squares of the deviations under the latter constraint, the following estimator is obtained in this case (OLS-estimator with constraints given by inequalities):

$$b_H = \begin{cases} b & \text{if } a_1 \leq b \leq a_2 \\ a_2 & \text{if } b > a_2 \\ a_1 & \text{if } b < a_1 \end{cases} \tag{32}$$

Denoting $C_1 = a_1 - \beta_1, C_2 = a_2 - \beta_1$ (with $C_1 \leq 0; C_2 \geq 0$) and using the results from the theory of moments of the functions of random variables, we will obtain the mean square of the distance from the estimate $b_H$ to $\beta_1$:

$$P^2(b_H) = 0.5 \tilde{\sigma}^2 \left[ \Phi(C_2 / \tilde{\sigma}\sqrt{2}) - \Phi(C_1 / \tilde{\sigma}\sqrt{2}) \right] +$$
$$+ 0.5 C_1^2 \left[ 1 + \Phi(C_1 / \tilde{\sigma}\sqrt{2}) \right] + 0.5 C_2^2 \left[ 1 - \Phi(C_2 / \tilde{\sigma}\sqrt{2}) \right] + \tag{33}$$
$$+ (C_1 \tilde{\sigma} / \sqrt{2\pi}) \exp(-C_1^2 / 2\tilde{\sigma}^{2)}) - (C_2 \tilde{\sigma}) / \sqrt{2\pi} \exp(-C_2^2 / 2\tilde{\sigma}^2),$$

where $\Phi(x) = (2/\sqrt{\pi})\int_0^x \exp(-t^2)\,dt$ is the Laplace integral. Let us construct now the resultant estimate $\overset{\bullet}{b}$, based on the use of the optimal quasi-estimator (17), for the case of additional information under consideration.

Denote by $b_1$ the minimal estimate from the two alternative ones, and by $b_2$ - the maximal one. Then the estimate $\overset{\bullet}{b}$ has the following form:

$$\overset{\bullet}{b} = \begin{cases} b_1 \text{ if } a_1 \leq b \leq a_2 \ \& \ a_1 < b_1 < a_2 \ \& \ b_2 > a_2 \ | \ b \geq a_2 \ \& \ a_1 < b_1 < a_2 \\ b_2 \text{ if } a_1 \leq b \leq a_2 \ \& \ b_1 < a_1 \ \& \ a_1 < b_2 < a_2 \ | \ b \leq a_1 \ \& \ a_1 < b_2 < a_2 \\ a_1 \text{ if } b < a_1 \ \& \ b_2 < a_1 \ | \ b < a_1 \ \& \ b_2 > a_2 \\ a_2 \text{ if } b > a_2 \ \& \ b_1 > a_2 \ | \ b > a_2 \ \& \ b_1 > a_2 \\ b \text{ if } a_1 \leq b, b_1, b_2 \leq a_2 \ | \ a_1 \leq b \leq a_2 \ \& \ b_1 < a_1 \ \& \ b_2 > a_2 \end{cases} \quad (34)$$

Where the sign $\&$ denotes the logical connective "and" and the sign $|$ – the logical "or."

For the estimate $\overset{\bullet}{b}$, an analytical expression for the mean sum of the error squares is also obtained for the case of the normal distribution of the error $\varepsilon$. Introduced as the sum of double integrals, it is reduced to the sum of one-dimensional integrals with the integrands containing incomplete gamma functions. This result is not given here because of its bulkiness. Besides, our task was studying cases with deviations from the starting assumptions. Therefore, in the sequel we will use the method of statistical modelling – which, with the modern capacities of the computing and of the software, is not less universal than computations based on theoretical relations.

Let us now consider the results of comparison of the OLS-estimator (32), under constrains given by the inequalities, – with the proposed method (34). The problem for $n = 6$, $\beta_1 = 1$ was studied. As before, the length of the application in this statistical testing was 10,000. The results of this statistical modelling are given in the tables 4-5.

It can be seen from the tables that, in all cases, (with the symmetric and asymmetric *a priori* region, relative to the parameter $\beta_1$) the proposed estimator has substantial advantage.

Let us also note that, for the canonical and non-canonical constrains of the form $\beta^T \beta = c^2$, or $\beta^T \beta \leq c^2$, which are used in the proof of the optimality of the ridge regression in [15], in our case, the problem of choosing the better of the two alternative estimates becomes trivial.

*Table 4. Comparative data for the mean squared distances*

*Normal distribution $\varepsilon \in N(0,1)$, for the OLS-estimator $\tilde{L}^2 = 0.1667$*

|  | Constraints | | | | |
|---|---|---|---|---|---|
| *Method* | $0.3 \leq \beta_1 \leq 1.7$ | $0.5 \leq \beta_1 \leq 1.5$ | $0.8 \leq \beta_1 \leq 1.2$ | $0.6 \leq \beta_1 \leq 1.7$ | $0.3 \leq \beta_1 \leq 1.4$ |
| *OLS-estimator with constraints* | 0.14140 | 0.10861 | 0.02985 | 0.11375 | 0.11157 |
| *Proposed method* | 0.07649 | 0.05126 | 0.01925 | 0.06337 | 0.06421 |

*Table 5. Comparative data for the mean squared distances*

*Uniform distribution $\varepsilon$ is in the interval ( -2,2) , for the OLS-estimator $\tilde{L}^2$ = 0.2216*

|  | Constraints | | | | |
|---|---|---|---|---|---|
| *Method* | $0.3 \leq \beta_1 \leq 1.7$ | $0.5 \leq \beta_1 \leq 1.5$ | $0.8 \leq \beta_1 \leq 1.2$ | $0.6 \leq \beta_1 \leq 1.7$ | $0.3 \leq \beta_1 \leq 1.4$ |
| *OLS-estimator with constraints* | 0.17698 | 0.12618 | 0.03148 | 0.13278 | 0.13259 |
| *Proposed method* | 0.08337 | 0.06216 | 0.02019 | 0.07532 | 0.07527 |

*Table 6. Comparative data for the mean squared distances*

*Mixed distributions $\varepsilon_1 \in N(0,1)$ – 80 % and $\varepsilon_2 \in N(0,10)$ - 20 %, for OLS-estimator $\tilde{L}^2$ =2.91*

|  | Constraints | | | | |
|---|---|---|---|---|---|
| *Method* | $0.3 \leq \beta_1 \leq 1.7$ | $0.5 \leq \beta_1 \leq 1.5$ | $0.8 \leq \beta_1 \leq 1.2$ | $0.6 \leq \beta_1 \leq 1.7$ | $0.3 \leq \beta_1 \leq 1.4$ |
| *OLS-estimator with constraints* | 0.38672 | 0.21148 | 0.03749 | 0.26001 | 0.26175 |
| *Proposed method* | 0.18776 | 0.12357 | 0.02957 | 0.14707 | 0.14458 |

## 4 Conclusion

The considered approach allows one to substantially increase the quality of parameters estimation in the orthogonal and non-orthogonal, one-dimensional and multi-dimensional regression problems, using minimal volume of additional information, which is practically always available to the researcher. The experience in application of the quasi-estimator demonstrates that the additional information used, in particular *a priori* one, can be either formalizable or non-formalizable (intuitive). Additionally, all available types of information can be used simultaneously, with the increase in the probability of the correct choosing of the better estimate out of two alternatives provided by the quasi-estimator.

It is also worth mentioning that the quasi-estimator is extremely simple to apply and requires only minimal modifications when working with any computer software for standard regression.